\documentclass[12pt]{article}
\usepackage{amssymb}
\def\om{\omega}
\def\sq{\subseteq}

\def\proof{\par\noindent Proof\par\noindent}
\def\poset{{\mathbb P}}

\def\forces{{| \kern -2pt \vdash}}

\def\qed{\par\noindent QED\par}

\def\cc{{\mathfrak c}}

\def\Tau{{\mathrm{T}}}
\def\rmor{\mbox{ or }}

\def\dual#1{\widetilde{#1}}
\def\rmand{\mbox{ and }}
\def\aron{Aronszajn }
\def\ss{\mathsf{S}}

\newtheorem{theorem}{Theorem}
\newtheorem{lemma}[theorem]{Lemma}

\newtheorem{question}[theorem]{Question}

\begin{document}

\begin{center}
{\large A Nonhereditary Borel-cover $\gamma$- set}
\end{center}

\begin{flushright}
Arnold W. Miller\footnote{
Thanks to the Fields Institute for Research in Mathematical Sciences
at the University of Toronto for their support during the time this paper
was written and to Juris Steprans who directed the special program in set
theory and analysis.
\par Mathematics Subject Classification 2000: 03E50; 03E17
}
\end{flushright}

\begin{center}
Abstract
\end{center}

\begin{quote}
In this paper we prove that if there is a Borel-cover $\gamma$-set of
cardinality the continuum, then there is one which is not hereditary. 
\end{quote}


In this paper we answer some of the questions raised by Bartoszynski and Tsaban
\cite{BT} concerning hereditary properties of sets defined by certain Borel
covering properties.

Define. An $\om$-cover of a set $X$ is a family of sets such that
every finite subset of $X$ is included in an element of the cover
but $X$ itself is not in the family.

Define. A $\gamma$-cover of a set $X$ is an infinite family of sets such
that every element of $X$ is in all but finitely many elements of the
family.

Define. A set $X$ is called a Borel-cover $\gamma$-set iff every countable
$\om$-cover of $X$ by Borel sets contains a $\gamma$-cover.
 
These concepts were introduced by Gerlits and Nagy \cite{GN} for open covers.
 
Being a Borel-cover $\gamma$-set is equivalent to saying that for any
$\om$-sequence of countable Borel $\om$-covers of $X$ we can choose one
element from each and get a $\gamma$-cover of $X$ - this is denoted   
$\ss_1({\cal B}_\Omega, {\cal B}_\Gamma)$.  The equivalence was proved by
Gerlitz and Nagy \cite{GN} for open covers but the proof works also for Borel
covers as was noted in Scheepers and Tsaban \cite{ST}: 

Let ${\cal B}_n$ be Borel $\om$-covers of $X$. Since 
 $$\{U\cap V: U\in{\cal U}, V\in{\cal V}\}$$ 
is an $\om$-cover if ${\cal U}$ and ${\cal V}$ are,  we may assume that
${\cal B}_{n+1}$ refines ${\cal B}_{n}$.  Let $x_n$ for $n<\om$ be distinct
elements of $X$ and let 
 $${\cal B}=\{A\setminus \{x_n\}: n<\om, A\in {\cal B}_{n}\}$$ 
It is easy to check that ${\cal B}$ is an $\om$-cover of $X$.  Now let
${\cal C}$ be a $\gamma$-subcover of ${\cal B}$. Note that for any fixed $n$ at
most finitely many of the elements of ${\cal C}$ can be of the form 
$A\setminus \{x_n\}$.  By refining $\cal C$ we may assume at most one thing is
taken from each ${\cal B}_n$ and since they are refining we can fatten up
${\cal C}$ to take exactly one element of each ${\cal B}_n$.

Define. A family of subsets of $X$,  $\;{\cal U}$ is a $\tau$ cover of
$X$ iff every element of $X$ is in infinitely many elements of ${\cal U}$ and 
for every $x,y\in X$ at least one of the sets  
  $$\{U\in{\cal U}: x\in U,\; y\notin U\} \rmor
  \{U\in{\cal U}: x\notin U,\; y\in U\}$$
is finite.

Clearly any $\gamma$-cover is a $\tau$-cover.  These covers were introduced
in Tsaban \cite{tsaban}.

\begin{theorem} 
Suppose there is a Borel-cover $\gamma$-set of size  the continuum.   Then
there is a Borel-cover $\gamma$-set $X$ and subset $Y$ of $X$ which is not a
Borel-cover $\gamma$-set.   In fact, there is an open $\om$-cover of $Y$
with no $\tau$-subcover.
\end{theorem}

\proof
For $X\subset P(\om)$ let 
 $$\dual{X}=\{\om\setminus a: a\in X\}$$
be the dual of $X$, i.e., the set of complements of elements of $X$.
Let $P\sq [\om]^\om$ be a perfect set of
independent subsets of $\om$.  This means that for every disjoint pair
$F_1,F_2$ of finite subsets of $P$ the set
$$\bigcap F_1\cap\bigcap\dual{F_2}\mbox{ is infinite.}$$
Such a set was first constructed by 
Fichtenholtz, Kantorovich, and Hausdorff, see Kunen \cite{kunen}. 
To construct
one, let $Q=\{(n,s):n\in\om, s\sq P(n)\}$.  Define
$A_x=\{(n,s): x\cap n\in s\}$ for each $x\sq \om$ and
$$P=\{A_x:x\sq \om\}\sq P(Q)=2^Q.$$

Let $Z\sq P$ be a Borel-cover $\gamma$-set of
cardinality the continuum.

\bigskip
\noindent Claim. $Z\cup\dual{Z}$ is a Borel-cover $\gamma$-set.

proof: Let $\{B_n:n<\om\}$ be a Borel $\om$-cover of 
$Z\cup\dual{Z}$.  Then it is easy to see that
$\{(B_n\cap \dual{B_n}):n<\om\}$ is an $\om$-cover of
$Z$.  This is because if $(F\cup\dual{F})\sq B_n$ then
$F\sq (B_n\cap\dual{B_n})$.
 
Since $Z$ is a Borel-cover $\gamma$-set
there exists an $a\in [\om]^\om$ such that
$$\{(B_n\cap \dual{B_n}):n\in a\}$$ is a $\gamma$-cover
of $Z$. But then it is also a $\gamma$ cover of $\dual{Z}$.
This proves the Claim.

\bigskip
Let $X=Z\cup\dual{Z}$ and to pick $Y\sq X$ as required 
we will choose for
each $a\in Z$ to put either $a\in Y$ or $(\om\setminus a)\in Y$
(but not both).
Since $Z$ was a subset of $P$ and $P$ was independent we will have
that the intersection of any finite subset of $Y$ is infinite.  In
particular,  
$$ {\cal U}=\{U_n:n\in\om\}\mbox{ where }  
U_n=\{a\sq \om:n\in a\} $$
is an $\om$-cover of $Y$. But $\{U_n:n\in b\}$ is a
$\gamma$-cover of $Y$ iff $b\sq^* a$ for every $a\in Y$.
But this is easy to defeat.  Using that $Z$ has cardinality the
continuum let $Z=\{a_\alpha:\alpha<\cc\}$ and let
$[\om]^\om=\{b_\alpha:\alpha<\cc\}$. For each $\alpha$
if $b_\alpha\sq^* a_\alpha$ put $(\om\setminus a_\alpha)$ into
$Y$ and otherwise put $a_\alpha$ into $Y$. 

To construct $Y$ so that ${\cal U}$ has no  $\tau$-subcovers
can be done
by using two elements $a_0,a_1$ of $Z$  for each $b\in [\om]^\om$.
First note that the set $\{U_n:n\in b\}$ is a $\tau$-cover of $Y$ iff 
$b$ meets every element of $Y$ in an infinite set and for
every two elements $a_0,a_1$ of $Y$ either
$(a_0\cap b)\sq^* a_1$ or $(a_1\cap b)\sq^* a_0$.

Notation: $a^{(0)}=a$ and $a^{(1)}=\om\setminus a$. 

\bigskip\noindent
Claim. There exists $i,j$ in $\{0,1\}$ such that
\par (a) $b\cap a_i^{(j)}$ is finite or
\par (b) both $b\cap a_0^{(i)}\cap a_1^{(1-j)}$ and 
$b\cap a_0^{(1-i)}\cap a_1^{(j)}$ are infinite.

proof:  Assume case (a) fails for all $i,j$ in $\{0,1\}$.  
The four sets $a_0^{(i)}\cap a_1^{(j)}$ partition 
$\om$ into infinite sets since $a_0$ and $a_1$  are independent.  
If all four meet $b$ in an infinite set then we are done. So
assume that $b\cap a_0^{(i)}\cap a_1^{(j)}$ is finite for
some $i,j$. But
since $b\cap a_0^{(i)}$ is infinite it must be that
$b\cap a_0^{(i)}\cap a_1^{(1-j)}$ is infinite.  A similar
argument shows 
$b\cap a_0^{(1-i)}\cap a_1^{(j)}$ is infinite.
This proves the Claim.

\bigskip
To kill off the possibility of $b$ giving a $\tau$-subcover we
put $a_i^{(j)}$ into $Y$ in case (a) or put 
both $a_0^{(i)}, a_1^{(j)}$ into $Y$ in case (b).
This proves the Theorem.
\qed

\bigskip Remark. Tsaban points out the following corollary of our result. In
Problem 7.9 of Bukovsky, Reclaw, and Repicky \cite{brr} it is asked whether
every $\gamma$-set of reals which is also a $\sigma$-set is a hereditary
$\gamma$-set.  It is shown in  Scheepers and Tsaban \cite{ST} that every 
Borel-cover $\gamma$-set (more generally $\ss_1({\cal B}_\Gamma, {\cal
B}_\Gamma)$-set) is a $\sigma$-set.  Hence the answer to the problem is no.

\bigskip
The following result is due to Brendle \cite{brendle}.  Our proof
is a slight modification of a result of Todorcevic - see Theorem 4.1 of Galvin
and Miller \cite{GM} and is perhaps simpler.

\begin{theorem} \label{ch}(Brendle)
Assume CH. Then there exists a Borel-cover $\gamma$-set of size 
$\om_1$.
\end{theorem}

\proof 
The idea is to construct an \aron  tree 
$T\sq \om^{<\om_1}$ of perfect sets. 

We construct perfect subtrees $p_s\sq 2^{<\om}$ for
$s\in T$ and $X_s\sq [p_s]$ countable dense sets such that
\begin{enumerate}
\item if $s\sq t$ then $p_s \supseteq p_t$,
\item if $s$ and $t$ are incomparable then $[p_s]\cap [p_t]=\emptyset$,
\item if $\alpha<\beta<\om_1$ and $n<\om$ then for
every $s\in T_\alpha$ there exists $t\in T_\beta$ with $s\sq t$
and $p_s\cap 2^n=p_t\cap 2^n$, and
\item for every sequence $(B_n:n<\om)$ of Borel subsets of 
$2^\om$ there exists $\alpha<\om_1$ such that either
for some finite $F\sq \cup \{X_s:s\in T_{\leq \alpha+1}\}$
no $B_n$ covers $F$ or there exists an $a\in [\om]^\om$ 
such that $\{B_n:n\in a\}$ is a $\gamma$-cover of  
$$\cup \{[p_s]:s\in T_{\leq \alpha+1}\} 
\;\cup\; \bigcup \{X_s:s\in T_{\leq \alpha}\}$$
\end{enumerate}

After the construction is completed we will let $X=\cup\{X_s:s\in T\}$.
The last item guarantees that $X$ will be a Borel-cover $\gamma$-set.
The first three items are simply to guarantee that our construction
can continue at limit levels.  To do the last item we use the following
Lemma.

Define for any perfect tree $p$ and $s\in p$, 
$$p\langle s\rangle =\{t\in p: t\sq s \rmor s\sq t\}$$

\begin{lemma}
Suppose $\langle p_n:n<\om\rangle $ are perfect trees and 
$({\cal B}_n:n<\om)$ 
is a sequence of countable Borel $\om$-covers of 
$\cup_{n<\om}[p_n]$. Then 
there exists 
perfect pairwise disjoint subtrees $q_n\sq p_n$ and 
$\{B_n\in {\cal B}_n:n<\om\}$ which is a $\gamma$-cover 
of $\cup_{n<\om}[q_n]$.
\end{lemma}
\proof
We can begin by refining the $p_n$'s so that $[p_n]$ are pairwise disjoint.
So we may as well assume this to begin with.  
Also since Borel sets have the (relative) property of Baire with
respect to each perfect set, by passing to perfect subsets we may
assume that each of our Borel covers is an open cover.

Note that for  finite sequences $(k_i:i<n)$ and $(t_i\in p_{k_i}:i<n)$, there
exists a  $U_n\in {\cal B}_n$ and $r_i\supseteq t_i$ with $r_i\in p_{k_i}$ such
that
 $$\cup_{i<n}[p_{k_i}\langle r_i\rangle]\sq U_n$$
Using this observation it is easy to construct a fusion sequence which produces
the $q_n$ and the required $\gamma$-cover. This proves the Lemma. 

\qed

Let $T_\alpha=T\cap \om^\alpha$.

In our construction of the tree we start by assuming that
$\{{\cal B}_\alpha:\alpha<\om_1\}$ is a list containing all
countable families of Borel subsets of $2^\om$.  At limit ordinals
$\alpha$, we use the usual fusion arguments to produce $p_s$ 
for $s\in T_\alpha$.  We take care of condition
(4) as follows. 

Suppose by induction we have already constructed:
$(p_s :s\in T_{\alpha})$ and $(X_s:s\in T_{< \alpha})$.

To obtain condition (4) let
$\{x_n:n<\om\}=\cup\{X_s:s\in T_{<\alpha}\}$, and define
$${\cal B}^n_\alpha=\{B\in {\cal B}: \{x_i: i<n\}\sq B\}$$
 
If some ${\cal B}^n_\alpha$ is not an $\om$-cover of
$\cup\{[p_s]:s\in T_{\alpha}$, 
then there exists a finite subset of $\{x_n:n<\om\}$

which is not covered by any $B\in {\cal B}_\alpha$.  
In this case, we choose
$X_s\sq p_s$ so that this finite set is included
in $\cup\{X_s:s\in T_{\leq\alpha}\}$.  
We then choose $p_{sn}$ so
that $p_{sn}\cap 2^n=p_s\cap 2^n$, $p_{sn}\sq p_s$, and $[p_{sn}]$ 
for $n<\om$ pairwise disjoint. 
We don't need to worry about ${\cal B}_\alpha$ because it cannot
be an $\om$-cover of $X$.

So we may assume each ${\cal B}_\alpha^n$ is an $\om$-cover of 
$\cup\{[p_s]:s\in T_{\alpha}\}$.

Apply the Lemma to the sequence 

$$( p_s\langle t\rangle \;:\; 
s\in T_\alpha \rmand t\in  p_s )$$

Then for each $s\in T_\alpha$ and $n<\om$ let 
$$p_{sn}=\cup\{ q_{s,t}: t\in 2^n\cap p_s\}$$

In this case we can take each $X_{sn}\sq [p_{sn}]$ to be an 
arbitrary countable dense subset.  

\qed 
 
\bigskip

Theorem \ref{consistent} is probably known but we include its
proof here for completeness.

\begin{theorem}\label{consistent}
Suppose that $M$ is any countable standard model of ZFC.  Then
there exists a ccc poset $\poset$ in $M$ of size continuum such
if $G$ is any $\poset$-filter generic over $M$, then $X=M\cap 2^\om$
in $M$ is a Borel-cover $\gamma$-set in $M[G]$.  Note that forcing
with $\poset$ does not change the size of the continuum in $M[G]$.
\end{theorem}
\proof

This is really a corollary of result noted by Gerlitz and Nagy \cite{GN}
that assuming MA (or even just MA($\sigma$-centered)) that every
set $X$ of size less than continuum is a $\gamma$-set. 

Let $\{B_n:n<\om\}$ be an $\om$-cover of $X$. For each $x\in X$ let
 $$a_x=\{n:x\in B_n\}.$$  
The family $\{a_x:x\in X\}$ has the finite intersection
property.  So there is a well-known ccc poset of size $|X|$  (see Kunen
\cite{set}) which adds an infinite $a\in [\om]^\om$ such that 
$a\sq^* a_x$ for each  $x\in X$.  Then $\{B_n:n\in a\}$ is a
$\gamma$-cover of $X$.  To obtain the model $M[G]$ simply iterate continuum
many times, with the usual dovetailing argument to take care of all sequences
of Borel sets in $M[G]$. 

\qed

\begin{question}
Does MA imply there exists a Borel-cover $\gamma$-set
of size the continuum?
\end{question}

The theorems in this section show that it is consistent that the
classes $\ss_1({\cal B}_\Omega,{\cal B}_\Gamma)$,
$\ss_1({\cal B}_\Omega,{\cal B}_{\Tau})$, and
$\ss_{\rm fin}({\cal B}_\Omega,{\cal B}_{\Tau})$ are not hereditary.  

I don't know about the other classes in Bartoszynski and
Tsaban \cite{BT}, for example:

\begin{question}
Is the class $\ss_1({\cal B}_\Omega,{\cal B}_{\Omega})$ hereditary?
\end{question}

For the definitions of these classes see \cite{BT}.

\begin{flushleft}
Arnold W. Miller \\
miller@math.wisc.edu \\
http://www.math.wisc.edu/$\sim$miller\\
University of Wisconsin-Madison \\
Department of Mathematics, Van Vleck Hall \\
480 Lincoln Drive \\
Madison, Wisconsin 53706-1388 \\
\end{flushleft}

\end{document}